\documentclass[12pt,a4paper,doublesp]{amsart}
\usepackage[english]{babel}
\usepackage[all]{xy}

\advance\topmargin -1.5cm
\advance\oddsidemargin -1.5cm
\advance\evensidemargin -1.5cm

\newcommand{\be}[1]{\begin{equation} \label{eq:#1}}
\newcommand{\beq}{\begin{equation}}
\newcommand{\ee}{\end{equation}}

\newcommand{\bea}[1]{\begin{eqnarray} \label{eq:#1}}
\newcommand{\beaq}{\begin{eqnarray}}
\newcommand{\eea}{\end{eqnarray}}
\newcommand{\bean}{\begin{eqnarray*}}
\newcommand{\eean}{\end{eqnarray*}}

\newdimen\AAdi%
\newbox\AAbo%
\def\AArm{\fam0 \rm}%
\def\AAk#1#2{\setbox\AAbo=\hbox{#2}\AAdi=\wd\AAbo\kern#1\AAdi{}}%
\def\AAr#1#2#3{\setbox\AAbo=\hbox{#2}\AAdi=\ht\AAbo\raise#1\AAdi\hbox{#3}}%

\def\BBe{{\AArm I\!E}}%

\def\BBn{{\AArm I\!N}}%

\def\BBr{{\AArm I\!R}}%

\makeatletter

\@addtoreset{equation}{section}
\makeatother

\begin{document}


\newtheorem{thm}{Theorem}[section]

\newtheorem{lem}[thm]{Lemma}

\newtheorem{defn}[thm]{Definition}

\newtheorem{pro}[thm]{Proposition}

\newtheorem{cor}[thm]{Corollary}
\newtheorem{rem}[thm]{Remark}
%
%
%
\title{Penalized estimators for nonlinear inverse problems}

\author[]{Jean-Michel Loubes \and Carenne Lude\~na  }

\address{CNRS - Laboratoire de Math{\'e}matiques, UMR 8628\\
Universit{\'e} Paris-Sud\\
B{\^a}t 425, 91405 Orsay cedex France \vskip .1in}
\address{Departamento de Matem\'aticas, IVIC, Venezuela \vskip .1in}
\email{Jean-Michel.Loubes@math.u-psud.fr}

\email{cludena@ivic.ve}

\begin{abstract}
 In this article we tackle the problem of inverse nonlinear ill-posed
problems from a statistical point of view. We discuss the problem
of estimating the non observed function, without prior knowledge
of its regularity,
 based on noisy observations over
subspaces of increasing complexity.  We consider two estimators: a
model selection type estimator for
   the nested and the non nested cases, and a regularized adaptive estimator.
   In both cases, we prove consistency of the estimators and give their rate of convergence.
\end{abstract}



\maketitle

\section{Introduction}
The term inverse problem is   used to denote a wide area of
problems arising from both pure and applied mathematics. It deals with
situations where one must recover information about a quantity or
a phenomenon under study, from measurements which are indirect and
possibly noisy as well. Intuitively, an inverse problem
corresponds to deriving a cause, given some of its effects. This
kind of problems confront practical
scientists on a daily basis.\\

In this article we are interested in recovering an unobservable
signal $x_0$ based on observations  \be{model}
y(t_i)=F(x_0)(t_i)+\varepsilon_i, \ee where $F:X\to Y$ is a
nonlinear functional, with $X,\: Y$ Hilbert spaces and
$t_i,i=1,\ldots,n$ is a fixed observation scheme. $x_0:  \BBr \to
\BBr$ is the unknown function to be recovered from the data
$y(t_i),\: i=1,\dots,n$. The regularity condition over the
unknown parameter of interest is expressed through the assumption
$x_0 \in X$. We assume that the observations $y(t_i) \in \BBr$
and that the observation noise $\varepsilon_i $ are i.i.d.
realizations of a certain random variable $\varepsilon$.
Throughout   the paper, we shall denote ${\bf
y}=(y(t_i))_{i=1}^n$. We assume $F$ is Fr{\'e}chet differentiable
and ill posed in the sense that  our noise corrupted observations
might lead to large deviations when trying to estimate $x_0$. In
a deterministic framework, the statistical model \eqref{eq:model}
is formulated  as the problem of approximating the  solution of
$$F(x)=y,$$  when $y$ is not known, but is only available through  an approximation $y^\delta$, $$\|y-y^\delta\|\leq \delta.$$
It is important to remark that whereas in this case consistency of the estimators depends on the approximation parameter $\delta$,
in (\ref{eq:model})
it depends on the number of observations $n$.

In general, in the linear case, the best $L^2$ approximation of $x_0$ is $x^+=F^{+}y$,
where $F^{+}$ is the Moore-Penrose (generalized) inverse of $F$,
does not depend continuously on the left-hand side $y$. In
particular, for the ill-posed case,  unboundedness of $F^+$
  entails that $F^+(y^\delta )$ is not close to $x^+$.
 Hence, the inverse operator needs to be, in some sense, regularized. In the nonlinear case, under ill-posedness of the initial problem, we will always mean that the solutions do not depend continuously on the data, but it is no so clear to find a criterion to decide whether a nonlinear problem is ill-posed or not. Thus, if $F$ is a compact operator, local injectivity around $x^+$ a solution of $F(x^+)=y$, implies ill-posedness of the problem, as it is quoted in \cite{engl96}. Regularization methods replace an
ill-posed problem by a family of well-posed problems. Their
solution, called regularized solutions, are used as approximations
to the desired solution of the inverse problem. These methods
always involve some parameter measuring the closeness of the
regularized and the original (unregularized) inverse problem.
Rules (and algorithms) for the choice of these regularization
parameters as well as convergence properties of the regularized
solutions are central points in the theory of these methods, since
they allow to find the right balance between stability and
accuracy.\vskip .1in When $F$ is linear, the statistical problem
has been extensively studied, although in general efficient
parameter choice is still under active research. Two main types
of estimators have been considered. First regularized estimators
such as Tikhonov type estimators, then non linear thresholded
estimators. The first approach has been studied in great detail.
An interesting early survey of this topic is provided by
O'Sullivan in \cite{osulliv}. In this setting, the main issues
are what kind of regularizing functional should be considered and
closely related what the relative weight of the regularizing
functional should be. More recently, Mair and Ruymgaart in
\cite{Ruym96} studied different regularized inverse problems and
proved the optimality of the rate of convergence for their
estimators. Special attention has been devoted in this setting
when considering a Singular value decomposition (SVD) of operator
$F$. We cite  the recent work in this direction
  developed by Cavalier and Tsybakov in \cite{tsy} or Cavalier,
Golubev, Picard and Tsybakov in \cite{cgpt}. The second approach
has its most popular version in the {\em wavelet-vaguelet}
decomposition introduced by Donoho \cite{donho}. In this case the
main issue is finding an appropriate basis over which $F^+$, the
generalized inverse, is almost diagonal. This idea is further
developed by    Kalifa and Mallat \cite{mallat1} who introduce
{\em mirror wavelets}. Closely related, Cohen, Hoffmann and Reiss
in \cite{galerkin} construct an adaptive thresholded estimator
based on Galerkin's method.\\ \indent However,  scarce
statistical literature exists when $F$ is non linear. Among the
few  papers available, we point out the works  \cite{sulli90} or
\cite{snieder} where some rates are given. A different type of
approach is developed in Chow and Khasminskii \cite{chow} for
dynamical inverse problems.
 Moreover Bissantz et al. in \cite{munk} discuss in their work a nonlinear version of
the method of regularization (MOR). In practical situations, a
linearization of the inverse operator is often performed, though
the effects of this linearization are seldom studied. In our
paper, under certain assumptions controlling the nonlinearity  of
the operator, we are able to give precise rates of convergence.
\vskip .1in Yet, nonlinear inverse problems are very common in
practice, since they often arise when studying the solution of a
noisy differential equation. Optical or diffuse tomography refers
to the use of low-energy probes to obtain images of highly
scattering media. The inverse problem for one of the earliest
models of optical tomography amounts to reconstructing the
one-step transition probability matrix, defined by a system of
nonlinear equations, see for instance the work of Gr{\"u}nbaum in
\cite{3dtomo}. The issue of tomography is also often tackled in a
Bayesian framework (see for instance \cite{Kaipio}). However this
method requires a large number of regularized inversions of the
operator and a strong a priori over the data, while a frequentist
adaptive method can solve this issue in a more efficient way. In
electromagnetics, the observations also follow this framework,
\cite{electro}. Moreover, a large class of problems in economy are
related to price variations under constraints, which also can be
modelled using differential equations. It is also the case in
econometrics when studying statistical issues with heterogeneity.
Adding an instrumental variable is equivalent to transforming the
unknown direct problem into a
nonlinear inverse problem, see for instance the work of Darroles, Florens and Renault in \cite{florens}. \\
We must remark, however, that in the deterministic case, this
problem has been tackled. The authors use mainly $L^2$
regularized Tikhonov type estimators and show that they provide a
stable method for approximating the solutions of a nonlinear
ill-posed inverse problem. For general references, we refer to
the work of Engl, Hanke and Neubauer in \cite{engl96},
\cite{engl2000}, Tautenhahn in \cite{tauten1} or Tikhonov and
Leonov in \cite{tikko}. The choice of the regularization sequence
is here crucial. On the one hand, this issue is often practically
solved by numerical methodology which relies on grid methods and
recursive algorithms, see for instance \cite{bk}, \cite{gradient},
\cite{gradientb} or \cite{iter}. This point of view is close to
cross-validation methods in a statistical framework, see also
\cite{Ruymb}.
 On the other hand, a priori optimal choices of the smoothing sequence are given in the work of
\cite{chino}, leading the way to an adaptive estimation of $x$.
Such techniques,
 as well as model selection techniques can be used in a random framework.\vskip .1in
Our goal in this article is to estimate the parameter of interest
$x_0$, when $F$ is an ill-posed but known operator. We aim at
using complexity regularization methods to construct this
estimator. We will deal with a large class of operators, linear
but also non linear operators that still undergo some assumptions
that will be made precise later. Moreover, we want this estimator
to achieve optimal
 rates of convergence when the smoothness of the true solution is not known
a priori. From the statistical literature (\cite{wahba}, \cite{osulliv}) it
is clear that we have to ask for some kind of penalization in
order to obtain satisfying results when estimating the inverse problem.\vskip .1in Hence we consider penalized M-estimators
minimizing quantities of the form
\begin{equation} \label{nou} \hat{x}_n= {\rm arg} \min_{x \in {\mathcal{X}}} \left( \gamma(y-F(x)(t))+
\alpha_n {\rm pen}(x,{\mathcal{X}}) \right),
\end{equation}
where ${\mathcal
{X}}$ is a specific set, $\gamma(.)$ is a loss-function, ${\rm
pen}(.,.)$ is a penalty over $x$ and/or ${\mathcal{X}}$, and $\alpha_n$ is a
decreasing sequence all of which   will be defined precisely later. The idea of penalized M-estimators is to find
 an estimator close enough to the data, close in the sense defined by $\gamma$ and with a regularity
 property induced by the choice of the penalty ${\rm pen}$. The smoothing sequence $\alpha_n$ balances the two terms.
  The greater $\alpha_n$, the smoother  the estimator will be,
  while the smaller $\alpha_n$ the closer the estimator will be to the
  data, maybe leading to a too rough estimate.\vskip .1in
 Two main types of penalties are considered: dimension penalties that will lead to model selection type
 estimators  and regularity penalties which will lead to   regularized estimators.\begin{itemize}
\item In the first case the idea,  as
developed in \cite{barron},   consists in looking at a sequence of
 subspaces of increasing complexity (for now complexity will be defined in terms of dimension and
 we shall look at a sequence of linear subspaces to be defined later).
 This is usually considered a discretization method. The estimator is a projection estimator and the choice of the set
 is defined by the penalty. \item
In the second case,
   a regularization functional ${\rm pen}(x,\mathcal{X})=J(x)$ (typically a quadratic functional) is introduced and   appropriate weights are considered
   for this term. Actually, this method is known as Tikhonov regularization when the penalty is quadratic, see for instance \cite{tikko}.
 \end{itemize}

In both cases, the choice of either the collection of spaces or the smoothing sequence determines the behaviour of the estimator.
Since we want to construct adaptive estimators,
we investigate choices that do not depend on an a priori knowledge of the regularity of the true function $x_0$.
Indeed, if the unknown function $x_0$ is of known regularity it is quite
simple to estimate
 good discretization or regularization schemes a priori: find a subspace or find a
regularization parameter such  that the the error is smaller that
a certain prescribed threshold. However, usually this
  is not possible since the smoothness of the solution
is unknown, and adaptive methods must be used. Adaptivity means here, that the construction of
the estimator does not require knowing beforehand the regularity of the function of interest to be recovered $x_0$. In the inverse
problems literature this is known as a posteriori methods. But, we do assume that the inverse operator is known as well
as some assumptions, such as its degree of ill-posedness.\\
  In model (\ref{eq:model}) the assumptions entail that the variance grows with  the complexity:
  this is an important difference with standard numerical analysis.
  Indeed, this difference yields other optimal rates  which are usual in statistics (\cite{galerkin}, \cite{tsy})
  and which  we will give below.\vskip .1in
  Finally, what is the price to pay here for nonlinearity ? Since we will consider a linear expansion of the operator
  $F$ in a neighborhood of $x$, the introduction of non-linearity requires controlling the linear part of the Fr\'echet differential operator
  $F^{'}$ in balls around the true solution $x_0$. As  opposed to linear problems, this fact entails the need of
  finding a ''good'' initial guess which we  shall denote $x^\star$. Moreover, the ill posedness of the problem
  requires relating the non linearity to the smoothing  properties of $F^{'}(x_0)$, see condition {\bf AF} below.
  We remark that this kind of condition is at the heart of probabilistic control of noise amplification. We show,
  that, under such restrictions, the nonlinearity of the problem does not change the rate of convergence and we still
  are able to build adaptive estimators.
\vskip .1in The article is divided into five main parts. In
Section \ref{sintro} we describe our general framework.   In
section \ref{smodel} we discuss discretization methods
(projection estimators) and find optimal rates for ordered and
non ordered selection.   In section \ref{stiko} we tackle
regularization methods: we prove optimality of an adaptive
Tikhonov like estimator.  Section \ref{slemmas} is devoted to
technical lemmas which will be useful in the paper.

\section{Presentation of the problem}  \label{sintro} \vskip .1in

In this section we introduce general notation and assumptions.
These include   standard concentration assumptions over the
observation noise and some restrictions over the class of
operators $F(.)$.
\subsection{General assumptions
} $ $ \vskip .1in
 Recall that we want to estimate a function $x_0: \BBr \to \BBr$. It is important to stress  that the observations depend on a
 fixed design $(t_1,\ldots,t_n) \in \BBr^n$. This will require
 introducing an empirical norm based on this design.
 Set $Q_n$ to be the empirical measure of the
covariables:
$$Q_n = \frac{1}{n} \sum_{i=1}^n \delta_{t_i} .$$ Here we have set $\delta$ the Dirac function.
The $L_2 ( Q_n )$-norm of a function $y \in Y$ is then given by
$$\| y \|_{n}= ( \int y^2 d Q_n )^{1/2} ,$$
and the empirical scalar product by
$$<y,\varepsilon >_n=\frac{1}{n} \sum_{i=1}^n \varepsilon_i y(t_i).$$
 Remark this empirical norm is defined over the observation space $Y$. Over
 the solution space $X$ we will consider the norm given by the
 Hilbert space structure. For the sake of simplicity, we will write $\|.\|_X=\|.\|$ when no confusion is possible  \\
 We   also  introduce certain standard assumptions on the observation
 noise
 \begin{description}
 \item[AN moment condition for the errors]$ $\\ $\varepsilon$ is a centered random variable satisfying the
 moment condition $\BBe (|\varepsilon|^p/\sigma^p)\le p!  /2$  and  $\BBe (\varepsilon^2)=\sigma^2$.
 \end{description}
\vskip .1in
 In the non linear case, our convergence analysis will ne a local one, hence it is   necessary to start with an
  initial guess of the solution. We require that this starting point $x^\star$ satisfies the following conditions:\begin{itemize}
\item This initial guess should  allow to construct a good
approximation of the unknown Fr\'echet derivative of the operator evaluated in $x_0$: $F'(x_0)$. For this, assume $F$ is Fr\'echet
differentiable and the range of $$DF(x_1,x_2)=\int_0^1
F'(x_1\theta + x_2(1-\theta))d\theta$$ remains unchanged in a neighbourhood of the guess solution $x^\star$, i.e in the ball $B_\rho(x^\star)=\{x, \: \|x-x^\star\| \leq \rho\}$, for a certain $\rho>0$. More precisely, we
 assume using the same assumption as in \cite{Kaltenbacher} \vskip .1in
 \begin{description}
\item[AF control over the non linear part of the differential operator]$ $\\
  There are
$c_T$, a fixed linear operator $T$ (generally $T=F'(x^\star)$) and a linear operator depending on $x$ and $x'$, written $R(x,x')$ such that for $x, x'\in B_\rho(x^\star)$
$$F(x)-F(x')=TR(x,x')(x-x'),$$
with $\|I-R\|\leq c_{T}$.
\end{description}
Hence, $T$ is a known bounded linear operator that can be seen as
some approximation to $F^{'}(x)$ in a neighbourhood
$B_\rho(x^\star)$, which we must be able to use in our
computations. Note also that, in contrast to that, provided a
bound holds for $\|I-R\|$, these linear operator needs not be
known explicitly. In the Section \ref{whynlin}, we discuss such
drastic restrictions for the operators and provide examples
satisfying to these assumptions.
\item Assume also, that the image by the linear operator $T$ of $x^\star$ and $x_0$ are close, in the sense that
\begin{description}
\item[IG identifiability condition]$ $\\
$ x_0-x^\star\in Ker(T)^{\bot}$,
 \end{description}
where $Ker(T))^{\bot}$ is the
orthogonal complement of the null space of the operator $T$. \\
The following approximation result   (\cite{Kaltenbacher}) assures uniqueness
of the sought solution if the initial guess is sufficiently close.
  \begin{lem} Assume {\bf AF} holds with $c_T<1/2$ and assume
    for $x, x^\prime \in B_\rho(x^\star)$, $F(x)=F(x^\prime)$ and
    $x-x^\prime \in Ker (T)^{\bot}.$ Then $x=x^\prime$.
    \end{lem}
This lemma guarantees the identifiability of the estimation
problem \eqref{eq:model} since the solution is uniquely chosen.
\end{itemize}
If $x^\star$ is such that $x_0\in B_\rho(x^\star)$, the local
behaviour of $F$ will be defined by operator $T$. We assume the regularity of the problem is defined by that of $F'(x_0)$ hence its approximation $T$:
 this linear operator acts  with a degree of ill-posedness defined by an index $p$, which comes from the fact
 that the operator is not compact and therefore its inverse is not $L^2$ bounded. This is generally expressed
 by the fact that $T$ maps $L^2$ into some Sobolev space $H_p$. This condition is quite natural when studying the ill-posedness of operators, see for instance \cite{galerkin}. That is the reason why, we assume $T$ acts along a Hilbert
 scale $H_s$.
\begin{description}
\item[IP ill posedness of the operator]$ $\\
There exists $p>0$ such that $F'(x_0)(H_s)=H_{s+p}$.
\end{description}
This property can be expressed in an equivalent way by the ellipticity property
 \begin{equation}
   \label{ellip}
   <Tx,x> \sim \|x\|^2_{H^{-p/2}}
 \end{equation}
where $H^{-p/2}$ stands for the dual space of the Sobolev space $H_{p/2}$.

 \subsection{Approximating subspaces}$ $ \vskip .1in

Estimating over all $X$ is in general not possible. We shall thus
assume we are
 equipped with a sequence $Y_1\subset Y_2\ldots \subset Y_m \ldots \subset Y$ of nested linear
 subspaces whose union is dense in $Y$.
We assume $$dim(Y_m)=d_m.$$

Denote the projection of $W$ over
 any subspace $Z$ by $\Pi_Z W$. Let $\Pi^n_{Y_m}$ stand for the projection in the empirical norm.
  Set $T_m= \Pi^n_{Y_m} T$, for the linear operator $T$ defined in {\bf
  AF}. Indeed we get the following diagram
$$\xymatrix{ X \ar[r]^T \ar[dr]^{T_m} & Y \ar[d]^{\Pi_{Y_m^n}} \\
& X_m}$$

As we are assuming that   over $Y$ we consider the
  empirical norm, and over $X$ the usual $L^2$ norm, the adjoint operator of $T_m$ with respect to such
  topology,
  $T^\star_m$, actually depends on the observation sequence $t_i$.
  However, we will usually drop this fact from the notation.

 For
 example, if $Y_m$ is generated by some orthonormal basis $\phi=(\phi_1,\dots,\phi_{d_m})$,
 with respect to the $L^2$ norm over $Y$, and $T={\rm Id}$, then
 $$\Pi_{Y_m}^n y=\sum_{j=1}^{d_m} y_{j,n} \phi_j ,$$ where $y_{j,n}=<\Pi_{Y_m}^n y,\phi_j>_n$
 are the solution to the projection problem under the empirical
 measure $Q_n$. Set $G=[\phi_j(t_i)]_{i,j}$ to be the Gram matrix associated to basis $\left(\phi_j \right)_{j\ge 1}$.
 Thus, $$y_{j,n}= (G_m^tG_m)^{-1}G_m^t (y(t_1),\ldots, y(t_n)).$$

  Define $X_m=T_m^* Y$.  Let $A^+$ stand for
the generalized inverse of a closed range operator $A$. Then, by
construction
$$\Pi_{X_m}=(\Pi^n_{Y_m}T)^+ \Pi^n_{Y_m} T.$$

Our goal is to estimate the unknown $x_0 $ by $x_m \in x^\star +
X_m$ in such a way that $F(x_m)$ is close to the observed $y$. By
assumption {\bf IG} this is saying we want to approximate
$Ker(T)^{\bot}$ by means of the collection $X_m$ in some sense
that is  adjusted  to the observation error.

Define
$$\nu_m=\|(\Pi^n_{Y_m}T)^+ \Pi^n_{Y_m}\|.$$
This quantity controls the amplification of the observation error
over the solution space $X_m$.   We have, \cite{bk},
   $$\nu_m\ge \gamma_m : =\inf_{v\in Y_m, \|v\|=1}\|T^* v\|.$$
Parameter $\gamma_m$ expresses the effect of operator $T^*$   over the
approximating subspace $Y_m$. In the case of $T$ acting over a
Hilbert scale $H_t$, we may assume $T(H_s)=H_{s+p}$ and
\begin{equation} \label{ladim}
\gamma_m=d_m^{-p} .
\end{equation}
On the other hand this term is related to the
goodness of the approximation scheme. Indeed, \cite{bk} it can be
seen that
$$\gamma_{m+1}\le \|(T^* (I-\Pi^n_{Y_m}))\|\le C \| I-\Pi^n_{Y_m}\|_n: = C \gamma^m.$$
as $\|T^*\|$ is bounded. Here $\gamma^m$ describes    the
  approximation properties of the sequence $Y_m$.   The next
  assumption requires that $\gamma_m\sim \gamma^m$, that is, that there exist
  constants $c_1, c_2$ such that $$c_2\le \gamma_m/ \gamma^m\le c_1.$$

\begin{description}
   \item[AS amplification error]$ $\\   Assume there exists a positive  constant $U$
   such that $\|I-\Pi_{Y_m}\|_n\le \sqrt{U}
   \gamma_m$.
   \end{description}
 \begin{rem} Assumption {\bf AS} thus establishes that the worst amplification
of the error over $X_m$ is roughly equivalent to the best
approximation over $Y_m$. If we think of $Y$ as smoothed by the
action of
 operator $T$, the above assumption is quite natural. In the case
 of an operator acting over a Hilbert scale, we may assume $Y=H_p$ and then
 $\gamma^m=d_m^{-p}$.
 \end{rem}

\subsection{Rates of convergence for inverse problems} $ $ \vskip .1in
 Consider the projection estimator $x_m$ over $X_m\cap B_\rho(x^\star)$. That is, $x_m$ is chosen in a
 such a way as to minimize
 $\|y-F(x_m)\|_n$. The estimation
error can be thus bounded by \be{opt1}\|x_0-x_m\|\le
\|(I-\Pi_{X_m})x_0\|+ \frac{\|\Pi_{Y_m}\epsilon\|}{\nu_m} .\ee

Since $Y_m$ is a sequence of linear subspaces $\BBe
\|\Pi_{Y_m}\epsilon\|_n^2=O(d_m/n)$. Thus rates are of order
$$\|(I-\Pi_{X_m})x_0\|+\frac{d_m^{1/2}}{n^{1/2}\gamma_m}.$$ This rate
depends on the ill-posedness of the operator and the approximation
properties of $X_m$. In some cases these are known precisely and
in others they can be deduced from the properties of the solution
$x_0$. One such case is the following source assumption
encountered typically in the inverse problems literature

\begin{description}
\item[SC source condition]$ $\\ There exists $0<\nu\le 1/2$ such that $x_0\in Range((T^*T)^\nu)=\mathcal{R}((T^*T)^\nu)$
\end{description}
Indeed consider
$$A_{\nu,\rho}=\{ x \in X, x=
(T^*T)^{\nu} \omega, \|\omega\| \leq \rho\}$$ where
$0\leq\nu\leq\nu_0$, $\nu_0 > 0$ and use the further notation

\be{3} A_{\nu}=\bigcup_{\rho>0}
A_{\nu,\rho}=\mathcal{R}((T^*T)^{\nu})\ee \vskip .1in

These sets are usually called source sets, $x \in A_{\nu, \rho}$
is said to have a source representation. The requirement for an
element to be in $X_{\mu, \rho}$ can be considered as an
smoothness condition.\\
Then, following \cite{bk} and under {\bf SC} we have if $\nu\le 1/2$,
$$\|(I-\Pi_{X_m}) x_0\|\le \|I-\Pi_{Y_m}\|_n^{2\nu}=O(d_m^{-2\nu p}),$$
for a certain $p$.

This leads to the rate $$ \|x_m-x_0\|=O(n^{-\frac{2\nu p}{4\nu p+2p +1}}).$$
Interpreting  this rate in the statistical literature reads
$s=2\nu p$: the regularity depends on the ill-posedness of the
problem.   In
the ill posed literature the error is not related to the
underlying dimension so that rates are different. Typically in a
Hilbert scale setting, if the true solution $x_0\in H_s$, optimal
rates are of order $O(n^{-s/(2s+2p+1)})$, see for example \cite{tsy}.

\subsection{Comments on the assumptions} \label{whynlin} $ $ \vskip .1in
We impose severe conditions [IP] and [AF] for the inverse operator.\vskip .1in First of all, we note that this paper handles the case of linear operators since the assumptions are fulfilled with $R={\rm Id}$ and $x^\star=0$. All the results are valid without any further assumptions.\vskip .1in
 For the non linear case, we cannot expect a general rate result without such a condition, which  measures the difference
  between the operator and its linear counterpart, i.e when $T=F^{'}(x^*)$, in terms related  to the ill posedness of the operator.
  \\The nature of the condition we impose   comes from the fact that the rates of convergence are drawn from the source conditions [SC].
   Hence, in the non linear case, we need to impose as quoted in \cite{engl2000} that $(F^{'}(x^*))^*(y-F(x^*))
    \in \mathcal{R}(F^{'}(x_0)^*).$ A necessary condition is given by $$ \mathcal{R}(F^{'}(x^*))^*)
    \subset \mathcal{R}(F^{'}(x_0)^*).$$ This condition assures we can find a linear operator $R(x,x^*)$, with $\|{\rm Id}-R\| \leq c_T$ for some
    constant $c_T$, such that for all $x \in B_\rho(x^*)$ we have
$$ F^{'}(x)=F^{'}(x_0) R(x,x^*).$$
   Other more standard assumptions  in the statistics literature for inverse problems,
   such as $F^{'}(.)$ being Lipschitz,  are not well suited to deal with   ill-posed problems. \\
As a matter of fact, as stated in \cite{hanke95} or \cite{Kaltenbacher}, the usual assumptions such as Lipsischitz condition (with constant $L$) on its Fr\'echet derivative
$$ \|F(x)-F(x_0)-F^{'}(x)(x-x_0)\| \leq \frac{L}{2} \|x-x_0\|^2, $$ implying an estimate of the first-order Taylor remainder, are not appropriate in this case. That is the reason why required range invariance conditions of the type
$$F(x)-F(x')=TR(x,x')(x-x'),\: \forall x,x'\in B_\rho(x^\star) $$
$$ \|R-{\rm Id} \|\leq c.$$ This means, as quoted before, that the range of the divided difference operators $DF(x,x')$ in
$$ F(x)-F(x')=DF(x,x')(x-x') $$ remains unchanged in $B_\rho(x^\star).$ \\
Such condition has been verified for many applications, see for
instance \cite{Englmedic}, \cite{hanke95}, \cite{K98} or
\cite{Deuf}. The alternative assumption, often used also in the
literature for numerical solution of non linear inverse problems,
is called the tangential cone condition
$$ \|F(x)-F(x^{'})-T(x-x^{'})\| \leq c \|T(x-x^{'}) \|, \: \forall x,x'\in B_\rho(x^\star).$$
Previous authors point out that assumption we have chosen in this paper, plays a role in certain parameter identification problems from boundary measurements where the tangential cone condition might be hard or impossible to verify.\\
\indent    As an example  consider the non linear Hammerstein inverse problem. The operator $F$ is given by:
\begin{align*}
 F: H_1[0,1] & \rightarrow L^2[0,1] \\
 x &  \rightarrow \int_0^t \phi(x(s))ds.
\end{align*}
If  $\phi$ is assumed to belong to $C^{2,1}(I)$ for all intervals $I
\subset \BBr,$ Hanke et al. in \cite{hanke95} prove that such
inverse problem is ill posed and that the operator undergoes the
previous restrictions.  In \cite{Kaltenbacher}, the assumptions are satisfied by the inverse groundwater filtration problem of identifying the transmissivity $a$ in
\begin{align} \label{theex}
  -\nabla (a \nabla u) & =f \quad {\rm in} \:\: \Omega \\
\nonumber  u & =g \quad {\rm on} \partial \: \: \Omega,
\end{align}
on a $C^2$ domain $\Omega \subset \mathbb{R}^3$, from measurements of the piezometric head $u$, where is the parameter-to-solution map $F:a \to u.$\vskip .1in
The approximation properties are quite natural having in mind Sobolev spaces of different order and piecewise polynomial approximations, as it is highlighted in \cite{Kaltenbacher}. There are some examples undergoing such conditions. For instance, in this previous work, it is also proved that the problem \eqref{theex} fulfills these assumptions with $p=2$, since the linearization $T$ of the operator $F$ acts in as smoothing a way as integrating twice.
\section{Complexity regularization}\label{smodel}
\subsection{Ordered selection}$ $ \vskip .1in
Consider to begin with that  $Y_m, m\in M_n$  is  a sequence of
nested subspaces. We define ordered selection as the problem of
choosing the best $m$ based on the observations. For this we will
construct {\em penalized} estimators that require finding the
first $m$ that minimizes $$\|\Pi_{Y_m}(y-F(x_m))\|^2_n + {\rm
pen}(m),$$ where ${\rm pen}(m)$ is an increasing function. From a
deterministic point of view this is essentially equivalent to
choosing $m$ based on the {\em discrepancy principle} (see
\cite{bk}for an application of the discrepancy principle to non
linear problems), however the fact that the error does not have a
finite energy and that the goodness of fit depends on the number
of observations introduces important changes both in the methods
of proof as in the definition of the estimator.

More precisely, define the estimator

\be{estip1} \hat{x}_{\hat{m}}=\arg\min_{m\in M_n} \arg\min_{x\in
x^\star+ X_m, x\in B_\rho(x^\star)}\| \Pi_{Y_m}^n(y-F(x))\|_n^2+{\rm pen}(m)
\ee where ${\rm pen}(m)>r (1+L)\sigma^2  d_m/n$, $r=2+\theta$, for some
$\theta>0$, and $L>0$.

Numerically, minimization in the above expression is more
complicated than it would be in the linear case because we must
calculate the projection matrix at each step. However, choosing
an efficient sampling scheme will do the job.

The next theorem says the above estimator is also efficient in
terms of the rates in equation (\ref{eq:opt1}) except for a
constant. As a matter of fact, the model selection estimator has a rate of convergence less or equal than the best rate achieved by the best estimator for a selected model.
\begin{thm}  There exist  constants $C(r,\sigma),  $ and $k(r,\sigma)$ such that with
probability greater than $1-2e^{- k u^{1/(2(p+1))}}$ \be{estip2}
 \|\hat{x}_{\hat{m}}-x_0\|^2\leq C(r,\sigma)  \inf_{m\in M_n}
(\|(I-\Pi_{X_m})x_0\|^2+\frac{{\rm pen}(m)}{\gamma^2_m
n})+\frac{u}{n} \ee
\end{thm}

Proof: For each $m$, assume $x_m=x^\star+z_m$, $z_m\in X_m$, and
$x_m\in B_\rho(x^\star)$. Set $$w( x_m)=\frac{\Pi_{Y_m
}^n(F(x_{m })-F(x_0))}{\|\Pi_{Y_m}^n(F(x_{m})-F(x_0))\|_n}.$$ We divide
the proof in a series of steps.
\begin{itemize}
\item Control of $d_{\hat{m}}$.
Recall the penalization is defined by ${\rm pen}(m)=
r(1+L)\sigma^2[d_m+1]/n$, with $2<r $ a certain constant.

 Following standard arguments we have

\bean &&\|\Pi_{Y_{\hat{m}}}(F(\hat{x}_{\hat{m}})-F(x_0))\|_n^2+{\rm pen}(\hat{m})\\
&\leq& |\|\Pi_{Y_
{m}}(F(x_m)-F(x_0))\|_n^2+2<\Pi_{Y_{\hat{m}}}(F(\hat{x}_{\hat{m}})-F(x_0)),\varepsilon>_n\\
&&-2<\Pi_{Y_{ {m}}}(F( {x}_{
{m}})-F(x_0)),\varepsilon>_n-\|\Pi_{Y_{\hat{m}}}\varepsilon\|_n^2+\|\Pi_{Y_{m}}\varepsilon\|_n^2
+{\rm pen}(m). \eean

 Let $0<\kappa<1$. Since $2ab\leq \kappa
a^2+\frac{1}{\kappa}b^2$, for any $a,b$ we have for any $   m$
and $x_m\in X_m$
\bean
&&2<\Pi_{Y_{m}}(F(x_m)-F(x_0)),\varepsilon>_n\\
&\leq& \kappa \|\Pi_{Y_{m}}(F(x_m))-F(x_0))\|^2_n
+\frac{1}{\kappa}|<w(x_m)  ,\varepsilon>_n|^2 .\eean

Set  for $x>0$, $t(m)=c[d_m+1]+(1+e)^{-1}x$ and assume $\kappa$
and $g,c$ are chosen in such a way that
$\frac{1}{\kappa}((1+g)+(1+1/g)c)=c_1<r(1+L)$. Remark, $\kappa$
can be chosen very close to one.

Thus, \bean
&&(1-\kappa)\|\Pi_{Y_{\hat{m}}}(F(\hat{x}_{\hat{m}})-F(x_0))\|_n^2+
(r(L+1)-c_1)
\frac{\sigma^2 [d_{\hat{m}}+1]}{n}\\
&\leq& (1+\kappa)\|\Pi_{Y_m}^n(F(x_m)-F(x_0))\|_n^2 +
\frac{1}{\kappa}|<w(\hat{x}_{\hat{m}} ),
\varepsilon>_n|^2-c_1\frac{\sigma^2([d_{\hat{m}}  +1])}{n}\\
&\: &- \|\Pi_{Y_{\hat{m}}}\varepsilon\|_n^2+\|\Pi_{Y_{m}}\varepsilon\|_n^2\\
&+& \frac{1}{\kappa}|<w( {x}_{ {m}} ),
\varepsilon>_n|^2-c_1\frac{\sigma^2( [d_{m}+1]  )}{n}
+(r(L+1)+c_1)\frac{\sigma^2[d_m+1]}{n}. \eean

Now since $\kappa<1$ and $c_1=c_1(e)<r(L+1)$, we have for fixed
$e>0$,

\bean && \frac{\sigma^2(d_{\hat{m}}+1)}{n}\\ &\leq&
\frac{2}{r-c_1}(\|F(x_m)-F(x_0))\|_n^2+{\rm pen}(m))\\
&+& \frac{2}{\kappa(r-c_1)}\sup_{m,\|u_m\|_n=1}(|<u_m,
\varepsilon>_n|^2 -\frac{c_1}{2}\kappa (\frac{\sigma^2[d_m+1]}{n})), \eean

for $u_m\in S_m$ a countable dense subset of $X_m$. Note that using Lemma \ref{lema0}
$$\sup_{ \|u_m\|_n=1}\left(\frac{1}{\kappa}<u_m,
\varepsilon>_n\right)=c \|\Pi_{Y_m}^n\varepsilon\|_n.$$  We will note  for any square matrix $A$ $$\rho^2(A)=\max\mbox{eigenvalue}(A^t A).$$
Also note that,
$$\rho(\Pi_{Y_m}^{n\  t}\Pi_{Y_m}^n)=1\quad {\rm and} \quad  Tr(\Pi_{Y_m}^{n\  t}\Pi_{Y_m}^n)=d_m.$$

Hence, by Lemma \ref{lema1}, there are   constants $d$ and $c_2$
such that \bean &&P\left[\sup_{m }\sup_{\|u_m\|=1 }\left(|<u_m,
\varepsilon>_n|^2-c_1/2\kappa (\frac{\sigma^2[d_m+1]}{n})\right) >\frac{x}{n} \right]\\
&\leq &\sum_{m }P( \sup_{\|u_m\|_n=1 }|< u_m,
\varepsilon>_n|^2-c_1/2\kappa (\frac{\sigma^2[d_m+1]}{n})>\frac{x}{n})\\
&\leq&   \sum_{m } \exp\{-\sqrt{d (  x+ c_2 L [ d_m+1] )}\}\le C_2
e^{-\sqrt{d x/2}}, \eean setting $C_2=\sum_m e^{-\sqrt{c_2 L [
d_m+1]/2}}$.

So that with probability greater than $1-C_2 e^{-\sqrt{d x/2}} $
we have
\beaq \label{cotaproba1}&& \frac{\sigma^2 [d_{\hat{m}}+1]}{n}\\
\nonumber&\leq& \inf_{m}\inf_{x_m\in x^\star +z_m}
\frac{2}{r(1+L)-c_1}(\|F(x_m)-F(x_0))\|_n^2+{\rm pen}(m)) +\frac{ x}{n}.
\eea
On the other hand let, for any given $m$, $\tilde{x}_m$
stands for the "projection" of $x_0$ over $x^\star+X_m$, i.e.
$\tilde{x}_m=x^\star+z_m$ is such that

\begin{equation}\label{eq:proj}\Pi_{Y_m}^n
F(\tilde{x}_m)=\Pi_{Y_m}^nF(x_0).\end{equation}

Let $K_2=\frac{1+c_T}{1-c_T}$ for $c_T$ defined in {\bf AF}. We
have the following lemma, the proof of which is exactly as that
of Lemma 2 in \cite{Kaltenbacher}.

\begin{lem}\label{lema2}
Assume $K_2$ and   $x^\star$ are
such that $$K_2
\|[I-\Pi_{X_m}](x_0-x^\star)\|+\|\Pi_{X_m}(x_0-x^\star)\|\le
\rho.$$ Then there exists $\tilde{x}\in x^\star+X_m$, such that
(\ref{eq:proj}) is satisfied.

\end{lem}

Clearly, this solution satisfies,

$$\Pi_{Y_m}^n(F(\tilde{x}_m)-F(x_0)-T(\tilde{x}_m-x_0)+T(\tilde{x}_m-x_0))=0,$$
so that since $(\Pi_{Y_m}^n T)^{+}\Pi_{Y_m}^nT=\Pi_{X_m}$

\bean \Pi_{X_m}(\tilde{x}_m-x_0) & = & (\Pi_{Y_m}^nT)^{+}\Pi_{Y_m}^n
T(R(\tilde{x}_m,x_0)-I)
(\tilde{x}_m-x_0))\\
\\ &=&     \Pi_{X_m}^n (R-I)(\tilde{x}_m-x_0)
\eean

On the other hand,

$$(I-\Pi_{X_m})(\tilde{x}_m-x_0)=(I-\Pi_{X_m})(x^\star-x_0),$$

and therefore, by condition [AF]

$$\|\tilde{x}_m-x_0\|\leq \frac{1}{1-c_T}\|(I-\Pi_{X_m})(x^\star-x_0)\|$$

and

\bean  \|F(\tilde{x}_m)-F(x_0)\|& = &\|(I-\Pi_{Y_m}^n)(F(\tilde{x}_m)-F(x_0))\|\\
&\leq& \gamma^m (1+c_T) \|\tilde{x}_m-x_0\| \\ & \leq & \gamma^m
\frac{1+c_T}{1-c_T}\|(I-\Pi_{X_m})(x^\star-x_0)\|, \eean where
$\gamma^m=\sup_{x} \frac{\|(I-\Pi_{Y_m}^n)Tx\|_n}{\|x\|}$.

  Now let $d_{m_{opt}}$ be
such that
$$ d_{m_{opt}}= {\rm arg} \min_{d_m} \left[ (\gamma^m)^2(\frac{1+c_T}{1-c_T})^2\|(I-\Pi_{X_m})(x^\star-x_0)\|^2+{\rm pen}(m) \right].$$
 Since we are looking at ordered
selection,
$$g(m)=\gamma_m^2(\frac{1}{1-c_T})^2\|(I-\Pi_{X_m})(x^\star-x_0)\|^2$$
is a decreasing sequence, so that  the minimizer must be such that
$g(m)={\rm pen}(m)$. Hence we have
$$P\left[(d_{\hat{m}}+1-\frac{4r(1+L) }{r(1+L)-c_1}(d_{m_{opt}}+1))_+>u\right]\leq   C_2 e^{-\sqrt{d u/2}},$$
and $d_{\hat{m}}\leq \frac{4
r(1+L)}{r(1+L)-c_1}(d_{m_{opt}}+1)-1+u$ with probability greater
than $1- C_2 e^{-\sqrt{d u/2}}$.

\item Error bounds: Set $$\Delta=\|F(\hat{x}_{\hat{m}})-F(x_0))\|_n^2-
\inf_{m}\inf_{x_m}\frac{2}{1-\kappa}(\|\Pi^n_{Y_m}(F(x_m)-F(x_0))\|_n^2+{\rm pen}(m))_+ \: .$$
Lemma \ref{lema1} also yields $$P(\Delta>  x/n)\leq
 C_2 e^{-d x/2}.$$

 \item We are now able to prove optimal rates for our estimator.
 For this we need to bound $\|\Pi_{Y_m}^n(F(\hat{x}_{\hat{m}})-F(x_0))\|_n$ from
 below.

 Let $\Omega_{dim}(u)$ be the set such that $d_{\hat{m}}\leq
 \frac{4r(1+L)}{r(1+L)-c_1}[d_{m_{opt}}+1]+u-1$. Let $\Omega_{fit}(u)$ be the set where
 $\Delta<u$. In this section we assume we are always in $\Omega_{dim}(u)\cup
 \Omega_{fit}(u)$. We require the following Lemma

 \begin{lem} Let $x\in {\mathcal R}(T^\star T)+x^\star$. There
 exists a constant $C$ such that
  $$\|\Pi_{
 Y_m} (F(x)-F(x^\star+\Pi_{X_m} (x-x^\star)))\|\leq C(1+c_T)\gamma^m
  \|(I-\Pi_{X_m})(x-x^\star)\|.$$

 \end{lem}
\begin{proof}
  Let $x-x^\star=T^\star y$, with $y\in {\mathcal R}(T)$. By
 definition   $$(I-\Pi_{X_m})x=T^\star (I-\Pi_{Y_m}^n)y=T^\star w,$$
 with $w\in Y_m^{ort}$.

 \bean &&\|\Pi_{Y_m}^n[F(x)-F(x^\star+\Pi_{X_m} (x-x^\star))]\|\\ &\leq&
  \|\Pi_{Y_m}^n TR((x-x^\star)-(\Pi_{X_m} (x-x^\star)-x^\star)))\|\\
 &=& \|\Pi_{Y_m}^n T R (I-\Pi_{X_m})(x-x^\star)\|\\ &\leq &
 \sup_{w\in Y_m^{ort}, \|w\|=1}\|\Pi_{Y_m}^n T R T^\star w\|_n
 \|(I-\Pi_{X_m})(x-x^\star)\|.
\eean

 The first term in the latter is in turn bounded by
 \bean &&\|\Pi_{Y_m}^n T R\|\sup_{w\in Y, \|w\|=1} \|T^\star(I-\Pi_{Y_m}^n) w\|\\
 &\leq & C(1+c) \sup_{x\in ker(T)^{ort}, \|x\|=1}\|(I-\Pi_{Y_m}^n)Tx\|_n
 \eean
 since $T^\star$ is the adjoint operator of $T$.
\end{proof} \vskip .1in
With this Lemma we have, for $x_m=x^\star+z_m$, $z_m\in X_m$

 \bean &&
  \|\Pi_{Y_m}^n(F(x_m)-F(x^\star+\Pi_{X_m}(x_0-x^\star)))\|^2\\
 &\ge &\|\Pi_{Y_m}^nTR \Pi_{X_m}(z_m-(x_0-x^\star))\|^2-[C(1+c_T)\gamma^m ]^2
 \|(I-\Pi_{X_m})(x_0-x^\star)\|^2.
 \eean

 On the other hand, we have
 \bean
 & &\|\Pi_{Y_m}^nTR \Pi_{X_m}(z_m-(x_0-x^\star))\|^2\\ &\geq& \|z_m-(x_0-x^\star)\|^2 (\inf_{y\in
 Y_m}\frac{\|\Pi_{Y_m}^n T R T^\star y\|}{\|T^\star y\|})^2\\
 &\geq& \|z_m-(x_0-x^\star)\|^2(\inf_{y\in
 Y_m}\frac{<\Pi_{Y_m}^n T R T^\star y,y>_n}{\|T^\star y\| \, \|y\|_n})^2\\
 &\geq&\|z_m-(x_0-x^\star)\|^2(\inf_{y\in
 Y_m}\frac{< R T^\star y,T^\star y>_n}{\|T^\star y\| \,
 \|y\|_n})^2\\
 &\geq&\|z_m-(x_0-x^\star)\|^2(\inf_{y\in
 Y_m}\frac{\|T^\star y\|^2-<(I-R) T^\star y,T^\star y>}{\|T^\star y\| \,
 \|y\|_n})^2\\ &\geq& \|z_m-(x_0-x^\star)\|^2 (1-c_T)^2(\inf_{y\in
 Y_m}\frac{\|T^\star y\|}{\|y\|_n})^2\\ &=& \|z_m-(x_0-x^\star)\|^2
 (1-c_T)^2 \gamma_m^2
 \eean

 Thus, since $(I-\Pi_{X_m})(x_0-x_m)=(I-\Pi_{X_m})(x_0-x^\star)$
 under assumption {\bf AS},

 \be{cota45}\|x_m-x_0\|^2\leq \frac{\|F(x_m)-F(x_0)\|^2}{(1-c_T)^2\gamma_m^2}
  +(1+\frac{C^2(1+c_T)^2
  U}{(1-c_T)^2})\|(I-\Pi_{X_m})(x_0-x^\star)\|.
 \ee
Inequality (\ref{eq:cota45}) is true for whatever $x_m\in
x^\star+X_m$. Over $\Omega(u)$, we know $$d_{\hat{m}}\leq
\frac{4r(1+L)}{r(1+L)-c_1}[d_{m_{opt}}+1]+u-1.$$ We distinguish
then two cases according to whether $\hat{m}<m_{opt}$ or not.
\begin{itemize}
\item In
the first case it is clear that $x_{\hat{m}}\in
x^\star+X_{m_{opt}}$ and $m$ in (\ref{eq:cota45}) can be replaced
by $m_{opt}$.
\item In the second case we have:
$$\|(I-\Pi_{X_{\hat{m}}})(x_0-x^\star)\|\leq
\|(I-\Pi_{X_{m_{opt}}})(x_0-x^\star)\|.$$
\end{itemize}
In any case we have,
over $\Omega(u)\cap \Omega_{fit}(u)$
 \bean &&
 \|\hat{x}_{\hat{m}}-x_0\|^2\\
 &\leq& \|(I-\Pi_{X_{m_{opt}}})(x_0-x^\star)\|^2(1+\frac{C^2(1+c_T)^2 U
 }{(1-c_T)^2}+
 \frac{2}{(1-\kappa)(1-c_T)^2}
 \max(U,\frac{(\gamma^{m_{opt}})^2}{\gamma_{\hat{m}}^2}))\\ &+&
 \frac{2 ( {\rm pen}(m_{opt})+u/n)}{(1-\kappa)(1-c_T)^2\gamma_{\hat{m}}^2}
  \\
  &\leq&\|(I-\Pi_{X_{_{m_{opt}}}})(x_0-x^\star)\|^2(1+
 \frac{ C^2 (1-\kappa)(1+c_T)^2+2}{(1-\kappa)(1-c_T)^2})
   \\ &+&\frac{2}{(1-\kappa)(1-c_T)^2}\frac{{\rm pen}(m_{opt})+u/n}{ \gamma_{m_{opt}}^2
   }\frac{1}{ d_{m_{opt}}^p }\left(\frac{4r(1+L)}{r(1+L)-c_1}[d_{m_{opt}}+1]+u-1\right)^p\\
   &\le & C(r)\left( \|(I-\Pi_{X_{_{m_{opt}}}})(x_0-x^\star)\|^2+\frac{{\rm pen}(m_{opt} )}{
   \gamma_{m_{opt}}^2}\right) +K(r) \frac{u^{p+1}}{n},
 \eean
 for some appropriate constants $C(r,\sigma)$ and $K(r,\sigma)$.
 Thus, \bean &&P((\|\hat{x}_{\hat{m}}-x_0\|^2-C(r,\sigma)\{ \|(I-\Pi_{X_{_{m_{opt}}}})(x_0-x^\star)\|^2+\frac{{\rm pen}(m_{opt}
 )}{
   \gamma_{m_{opt}}^2}\})_+>\frac{K(r,\sigma) u^{p+1} }{n})\\
   &&\le 2e^{- \sqrt{d/2u}},\eean
 which ends the proof.

\end{itemize}

\subsection{Non ordered selection }$ $ \vskip .1in

Ordered selection has the advantage of working directly on the
observation space. It has the disadvantage that the expansion of
the solution $x_0$ over the resulting subspace $X_m$ might not be
efficient. This introduces the need for non ordered selection, or
equivalently, for threshold methods. The combination of both
ill-posedeness and non linearity yields this a difficult problem.
Indeed, the former yields that it is no longer possible to work on the
observation space as this would require simultaneous control of
$\gamma_m$ and $d_m$. Working on the solution space  requires
considering the inverse of a certain matrix. The goodness of fit
of the estimator is then defined by the trace and spectral radius
of this inverse matrix restricted to the sequence of subspaces, which in turn depends on the degree of nonlinearity of the
problem.

More precisely, let $m_0$ be such that
$$\|(I-\Pi_{X_{m_0}})(x_0-x^\star)\|\le
\inf_{m}[\|(I-\Pi_{X_{m }})(x_0-x^\star)\|+\sqrt{\frac{d_m}{n}}\frac{1}{\gamma_m}].$$

This quantity can be chosen so as not to depend on the unknown
regularity of the solution $x_0$. Under assumption {\bf SC} the
above inequality is satisfied if the dimension of the set is such that $$d_{m_0}^{2\nu p}\leq n^\frac{2 \nu
p}{4\nu p+2p +1}.$$ Thus it is enough to choose $m$ such that
$d_{m_0}\le n^{1/2p}.$ Analogous results are obtained in the case
of Hilbert scales (\cite{galerkin}). On the other hand, if $m_0$ is
estimated as is section 2, we have $\|(I-\Pi_{X_{m_0}})(x_0-x^\star)\|$
satisfies the optimal rates with high probability.

For this fixed $m_0$ set $A_{m_0}=T_{m_0}^{+}\Pi_{Y_{m_0}}^n$. Let
$\{Y_m\}_m \subset  Y_{m_0}$ be a collection of not necessarily nested
subspaces. We will use the notation $m\in m_0$ to resume the imbedding of such subsets. Our goal is to
find the best subspace along this collection using penalized
estimation. \\For fixed $x$, $D_m(x)=\Pi_{X_{m_0}}R(x)\Pi_{X_m}$ is
a linear operator, $D_m(x):X_m\to X_{m_0}$. Let $S_m$ be the
matrix whose entries are defined by
$$S_m(i,j)=\sup_{x\in B_\rho(x^\star) }|(A_{m_0}^t D_m(x))(i,j)|.$$

 Set $\rho_m=\rho(S_m^t S_m)$  and $t_m=Tr(S_m^t S_m)$. Let
 $R_m= t_m/\rho_m$. Remark that under our assumptions, namely that the basis is orthonormal for the fixed design,
   both $ n\rho_m$ and $ nt_m$  do not depend on $n$. Introduce $L_m$ a certain weight factor and in the
 notation of lemma \ref{lema1}, define  $\Sigma_i=\Sigma_i(m_0)$, $i=1,2$
 by
\be{cond1}  \Sigma_1=\sum_{m\in m_0}     e^{-\sqrt{
 d/2 r L_m (R_m+1)}} \ee
and for any $q \geq 2$ and a constant $C_q$ depending on $q$ \be{cond2}  \Sigma_2=\sum_{m\in m_0}  C_q
 (n \rho_m\sigma^2/d)^q [ (d/2 r L_m (R_m+1))^{q-1/2}+(d/2 r L_m (R_m+1))^{q-1}]   e^{-\sqrt{
 d/2 r L_m (R_m+1)}} \ee
 As before we will consider penalized estimation. The penalty
 term in this case will be set to
 $${\rm pen}(m)=r \sigma^2(1+L_m) [t_m+\rho_m],$$ with $r>2$. We now define the estimator by

\be{esti3}\hat{x}_{\hat{m}}=x^\star+{\rm arg}\min_{m \in m_0} {\rm arg}\min_{x_m\in X_m}\|A_{m_0} (y-F(x_m))\|^2+{\rm pen}(m).
\ee

Then, we have the following result

\begin{thm}  Assume {\bf IG, AF} and {\bf SC} are satisfied. Assume (\ref{eq:cond1})  holds true. Then,
 if $0<\kappa<1$,  with
probability greater than $1-\Sigma_1 e^{-\sqrt{d/(2n\rho_{m_0})
u}} $ \bean &&\|\hat{x}_{\hat{m}}-x_0\|^2\le
\|(I-\Pi_{X_{m_0}})(x_0-x^\star)\|^2\\ &&+\frac{2}{
(1-c_T-\kappa)}\inf_m\{{\rm arg}\min_{x_m\in
x^\star+X_m}(1+c_T)\|x_m-x_0\|^2+{\rm pen}(m)\}+ \sigma^2 u/n\eean

and for $q \geq 1$ \bean &&\BBe [\|\hat{x}_{\hat{m}}-x_0\|^2]^q \le [
\|(I-\Pi_{X_{m_0}})(x_0-x^\star)\|^{2}\\ &&+\frac{2}{
 1-c_T -\kappa}\inf_m\{{\rm arg}\min_{x_m\in x^\star+X_m}(1+c_T)\|
x_0-x_m\|^2+{\rm pen}(m)\} ]^q +\frac{ \Sigma_2}{n^q} \eean

\end{thm}

\begin{rem}
The above result depends on two factors: $t_m$ and $R_m$.
However, again for fixed $x$, we have
$$A_{m_0}DF(x)w=A_{m_0}TR(x)w=\Pi_{X_{m_0}}R(x)w,$$ so that bounds can
be obtained if we know $DF$. In order to understand the penalty
term, assume $F$ is linear, $F=T$, and let $b_j, \phi_j, \psi_j$ be its
Singular value decomposition (SVD), that is, $T\phi_j =b_j \psi_j$
and $T^* \psi_j=b_j \phi_j$. Let $Y_m$   be the linear space
generated by $(\phi_j)_{j\in m}$ for a collection of indices $m$.
 Assume the Gramm matrix $G $ associated to $ (\phi_j)_{j\in m_0}$
 satisfies $c\le \rho(G)\le C$. Then, ${\rm pen}(m)$ is roughly
 proportional to $1/n\sum_{j\in m}1/b_j^2$ and $R_m$ is proportional
 to $1/n\sup_{j\in m }1/b_j^2\times \sum_{j\in m}1/b_j^2$.  In the
 general case, for any collection $Y_m$, ${\rm pen}(m)$ will depend on
\begin{itemize}
 \item The Gramm matrix $G$,
\item how operator $T$ acts along $Y_{m_0}$,
 that is, $\inf_{v\in Y_{m_0}}\|T^*v\|/\|v\|$
\item and on
 the non linearity of operator $F$, that is, $|t_m-Tr(A_{m_0}^t \Pi_{X_m})|.$
\end{itemize}
In the complete diagonal form (${\rm pen}(m)\sim 1/n \sum_j
1/b_j^2$), the
 penalty entails a hard thresholding scheme: choose a coefficient
 if $(\hat{x}_{\hat{m}})_j\ge \sqrt{C/b_j^2 n}$ for some constant
 which depends on $m_0$.

\end{rem}

\begin{rem} In the linear case, $F=T$ and $x^\star=0$, the minimization problem in (\ref{eq:esti3}) can be simplified offering
an important insight. Indeed, note that  this problem is actually equivalent to minimizing
\bean \hat{x}_{\hat{m}}&=& {\rm arg}\min_m {\rm arg}\min_{x_m\in X_m}\{-2<A_{m_0}y,A_{m_0}Tx_m>+\|A_{m_0}Tx_m\|^2\}+{\rm pen}(m)\\
&=& {\rm arg}\min_m {\rm arg}\min_{x_m\in X_m}\{-2<\Pi_{x_{m}}A_{m_0}y,x_m >+\| x_m\|^2\}+{\rm pen}(m).
\eean
So that, for each $m$,
$$x_{m,j}=<A_{m_0}y,e_j>=<{\bf y}, A_{m_0}^te_j>, \, j=1,\ldots,m.$$
Thus, $m$ is selected by minimizing
$$-\sum_{j\in m}x_{m,j}^2+r\sigma^2(1+L_m)[\sum_{j\in m}\lambda_j+\sup_j \lambda_j],$$
for $\lambda_j$ the eigenvalues of $A_{m_0}^t \Pi_{X_m}$, which is a hard thresholding scheme. Note that the ill conditioned
$A_{m_0}$ need not be applied to the observation vector ${\bf y}$.

In the non linear case the problem is equivalent to minimizing
$${\rm arg}\min_m {\rm arg}\min_{x_m\in x^\star+X_m}\{-2<A_{m_0}(y-F(x^\star)),\Pi_{X_{m_0}}R(x_m,x^\star)(x_m-x^\star)>
+\| x_m-x^\star\|^2\}+{\rm pen}(m).$$
Set ${\bf F(x^\star)}=(F(x^\star)(t_i))_{i=1}^n$. Hence,
\begin{align*}
(x_m-x^\star)_j & =<A_{m_0}(y-F(x^\star)),\Pi_{X_{m_0}}R(x_m,x^\star)e_j> \\
& =<{\bf y}-{\bf F(x^\star)}, A_{m_0}^t\Pi_{X_{m_0}}R(x_m,x^\star)e_j>, \, j=1,\ldots,m.
\end{align*}
Then $m$ is chosen as above. However, in this case the  problem must be solved numerically which is troublesome as $A_{m_0}$
is a badly conditioned matrix.
\end{rem}

Proof: from the definition for any $m$ and $x_m$,

\bean && \|A_{m_0} (F(x_0)-F(\hat{x}_{\hat{m}}))\|^2\le \|A_{m_0}
(F(x_0)-F(x_m))\|^2+{\rm pen}(m) \\ &&+2<\varepsilon,A_{m_0}^*A_{m_0}(F(x_0)-F(\hat{x}_{\hat{m}}))>_n
+ 2<\varepsilon,A_{m_0}^*A_{m_0}(F(x_0)-F(
x_m))>_n-{\rm pen}(\hat{m}).\eean

We have $$A_{m_0}(F(x_1)-F(x_2))=\Pi_{m_0}R(x_1,x_2)(x_1-x_2).$$
Hence, the left hand side is bounded from below by
$(1-c_T)\|\Pi_{X_{m_0}}(\hat{x}_{ \hat{m}}-x_0)\|^2$  and
$$\|A_{m_0} (F(x_0)-F(x_m))\|^2\le (1+c_T)\|\Pi_{X_{m_0}}(x_m-x_0)\|^2.$$
  Thus, \bean
&&(1-c_T)\|\Pi_{X_{m_0}}(\hat{x}_{ \hat{m}}-x_0)\|^2\\&&
\le (1+c_T)\|\Pi_{X_{m_0}}(x_m-x_0)\|^2 + {\rm pen}(m)\\
&&+ 2<\varepsilon,A_{m_0}^*R(\hat{x}_{ \hat{m}},x_m)(\hat{x}_{
\hat{m}}-x_m )>_n -{\rm pen}(\hat{m}).\eean For any $m,m^\prime$ set
$\Pi_{m\backslash m^\prime}=\Pi_{X_m\backslash X_{m^\prime}}$ and
$\Pi_{m\cap m^\prime}=\Pi_{X_m\cap X_{m^\prime}}$. With this
notation
$$\|x_{m^\prime}-x_m\|^2=\|\Pi_{m\cap m^\prime}(x_{m^\prime}-x_m)\|^2+\|\Pi_{m\backslash m^\prime}(x_{m^\prime}-x_m)\|^2+
\|\Pi_{  m^\prime\backslash m}(x_{m^\prime}-x_m)\|^2,$$
 and \bean && |<\varepsilon,A_{m_0}^*R(\hat{x}_{
\hat{m}},x_m )(\hat{x}_{ \hat{m}}-x_m)>|\\&&=
|<\varepsilon,A_{m_0}^*R(\hat{x}_{ \hat{m}},x_m
)\Pi_{m\cap\hat{m}}(\hat{x}_{ \hat{m}}-x_m)>_n+
<\varepsilon,A_{m_0}^*R(\hat{x}_{ \hat{m}},x_m
)\Pi_{m\backslash\hat{m}}(\hat{x}_{ \hat{m}}-x_m)>_n\\
&&+<\varepsilon,A_{m_0}^*R(\hat{x}_{ \hat{m}},x_m
)\Pi_{\hat{m}\backslash m}(\hat{x}_{ \hat{m}}-x_m)>_n|\\
&&\le \kappa \|\hat{x}_{ \hat{m}}-x_m \|^2+2/\kappa \|\varepsilon
S_{\hat{m}}\|^2+1/\kappa\|\varepsilon S_{m}\|^2.\eean

 The first term in the latter is bounded by $$\kappa [\|\hat{x}_{ \hat{m}}-(\Pi_{X_{m_0}}x_0+x^\star) \|^2+
 \| x_m -(\Pi_{X_{m_0}}x_0+x^\star)\|^2].$$ The proof then follows  directly from lemma
\ref{lema1}.

\section{Regularization}\label{stiko}
Crucial questions in applying regularization methods are
convergence rates and how to choose regularization parameters to
obtain optimal convergence rates.

Yet another approach is to consider a big enough subspace
$Y_{m_0}$ and in order to deal with the ill posedness of $(T^*
\Pi_{Y_m}^n)^+$  use Tikhonov  regularization methods.

  As in the last section assume that $m_0$ is such that
$$\|(I-\Pi_{X_{m_0}})(x_0-x^\star)\|\le
\inf_{m}[\|(I-\Pi_{X_{m}})(x_0-x^\star)\|+\sqrt{\frac{d_m}{n}}\frac{1}{\gamma_m}].$$
Next consider for a given $k \in \mathcal{K}$ a sequence $\alpha_k(n) \to 0$ as $n \to \infty$ and define the
following penalized estimator:
\begin{equation}
  \label{Mest}
  \hat{x}_{\alpha_k(n)}=x^\star+{\rm arg}\min_{x \in x^\star+X_m} \left[ \|\Pi_{Y_{m_0}}(y-F(x)) \|_n^2+ \alpha_k(n) \|
    x-x^\star \|^2 \right]
\end{equation}
We point out that choosing the smoothing sequence $\alpha_k(n)$ is the key point since it balances the two terms: if $\alpha_k$ is big the solution will be smooth but will not, in general, comply to
the observations. On the other hand, if $\alpha_k$ is small, the
solution might be too close to the noisy observations to yield a
good approximation of $x_0$. From the theory of inverse problems, we know that it is possible to choose a regularization sequence achieving the optimal rate of convergence, but this choice depends on $\nu$, which characterizes in a way the regularity of the solution. \\
As for the projection problem we
would like too choose $\alpha_k(n)$, among all the $\alpha_k(n), \: k \in \mathcal{K}$ based on the data in such a way
that optimal rates are maintained. This choice must also not depend on a priori regularity assumptions.

In the linear case, if $F=T$ for each fixed $\alpha_k(n)$, the expression \eqref{Mest} can be written in the following way:
\begin{equation} \label{Mest2} \hat{x}_{\alpha_k(n)}= x^\star+ {\rm arg}\min_{x \in x^\star+X_m}\|(T^t\Pi_{Y_{m_0}}T+\alpha_k(n)
I_{m_0})^{-1}T^t\Pi_{Y_{m_0}}(y-T(x))\|^2.
\end{equation}

However, in the nonlinear case both estimators are not the same.
Although in practice the second one  is more complicated (the
matrix to inverse might be  big) in order to select $\alpha_k$,
we will choose this version of the estimator, with $T$ defined in
assumption {\bf AF}. With this notation set
$$R_{\alpha_k(n)}=(T^t\Pi_{Y_{m_0}}T+\alpha_k(n)I_{m_0})^{-1}T^t\Pi_{Y_{m_0}}.$$

Now set $\gamma(x,\alpha_k)=\|R_{\alpha_k}(y-F(x))\|^2$ and
$${\rm pen}(\alpha_k)= r\sigma^2(1+L_k)
[Tr(R_{\alpha_k}^tR_{\alpha_k})+\rho^2( R_{\alpha_k})],$$ with $r>2$.

We choose $\hat{x}_{\alpha_{\hat{k}}}$ such that
$$ \hat{x}_{\alpha_{\hat{k}}} =x^\star+ {\rm arg}\min_{k,x \in X_{m_0}} \left( \gamma(x,\alpha_k(n))+{\rm pen}(\alpha_k(n)) \right).$$

Let $x_{\alpha_k}= x^\star+R_{\alpha_k}T(x_0-x^\star)$.

Set
$$\Sigma(d)= \sum_k
2\left[\sqrt{\frac{d
Tr(R_{\alpha_k}^tR_{\alpha_k})}{\rho^2(R_{\alpha_k})}} + 1\right]
\frac{d}{\rho^2(n R_{\alpha_k})}e^{-\sqrt{d L_k
[Tr(R_{\alpha_k}^tR_{\alpha_k})+\rho^2(R_{\alpha_k})]/
\rho^2(R_{\alpha_k})}},$$ for $d $   as in lemma \ref{lema1}.

We have the following result,

\begin{thm} For any $x\in x^\star+X_m$ and any $k$ such that $d_{m_0}\ge
\alpha_k^{-1/(2p)}$,  the following inequality holds true
\be{desprim}\BBe
\|\Pi_{X_{m_0}}(\hat{x}_{\alpha_{\hat{k}}}-x_0)\|^2\le
\frac{1}{1-c_T}\inf_{k \in
\mathcal{K}}[C(1+c_T)\|\Pi_{X_{m_0}}(x_{\alpha_k}-x_0)\|^2+2{\rm
pen}(\alpha_k)] +\frac{\Sigma(d)}{(1-c_T)n}, \ee
\end{thm}
Hence, the estimator is optimal in the sense that the adaptive estimator achieves the best rate of convergence among all the regularized estimators.
\begin{rem} The condition $d_{m_0}\ge
\alpha_k^{-1/(2p)}$ follows from the requirement that
$$\frac{Tr(R_{\alpha_k}^tR_{\alpha_k})}{\rho^2(R_{\alpha_k})}
=O(\alpha_k^{-1/(2p)}).$$
\end{rem}

\begin{proof} There exists a constant $C$ such that for any $k$,
$$(1-c_T)\|\Pi_{X_{m_0}}(x_1-x_2)\|^2\le \|R_{\alpha_k}(F(x_1)-F(x_2))\|_n^2\le C (1+c_T)\|\Pi_{X_{m_0}}(x_1-x_2)\|^2.$$
Thus, for any $x\in x^\star+ X_{m_0}$,
 \bean && (1-c_T)\|\Pi_{X_{m_0}}(\hat{x}_{\alpha_{\hat{k}}}-x_0 )\|^2\le C(1+c_T)
 \|\Pi_{X_{m_0}}(x-x_0 )\|^2 \\
&+ & 2 {\rm pen}(\alpha_k)+2\sup_{k}[\|R_{\alpha_k}\varepsilon \|^2_n-{\rm pen}(\alpha_k)].\eean
 The result follows directly from Lemma \ref{lema1}.\\
In the linear case $F=T$, we get the following proof:
For any $x_{\alpha_k}$ and any $k \in \BBn$

\bean &\|R_{\alpha_{\hat{k}}}(y - T\hat{x}_{\alpha_{\hat{k}}})\|^2
+ pen(\alpha_{\hat{k}}) \leq \|R_{\alpha_k}(y - T
x_{\alpha_k})\|^2 + pen(\alpha_k) \eean and \bean &
\|R_{\alpha_k}(y - T x_{\alpha_k})\|^2= \|R_{\alpha_k}T({x}_0
- x_{\alpha_k})\|^2 + 2 \langle R_{\alpha_k} T ({x_0} -
x_{\alpha_k}), R_{\alpha_k}\varepsilon \rangle + \|R_{\alpha_k}
\varepsilon\|^2 \eean \vskip .1in

Thus, following standard arguments we have

\bean &&\|R_{\alpha_{\hat{k}}}T({x_0}-\hat{x}_{\alpha_{\hat{k}}})\|^2\\
&\leq& \|R_{\alpha_k}T({x_0}-x_{\alpha_k})\|^2 - 2<
R_{\alpha_{\hat{k}}}
T({x_0}-\hat{x}_{\alpha_{\hat{k}}}),R_{\alpha_{\hat{k}}}\varepsilon>\\
&&+2<R_{\alpha_{k}}
T({x_0}-x_{\alpha_{k}}),R_{\alpha_k}\varepsilon>-
\|R_{\alpha_{\hat{k}}}\varepsilon\|^2+
\|R_{\alpha_k}\varepsilon\|^2 +{\rm pen}(\alpha_k)+{\rm
pen}(\alpha_{\hat{k}}). \eean \vskip .1in

Let $0<\kappa<1$. Since $2ab\leq \kappa a^2+\frac{1}{\kappa}b^2$,
for any $a,b$ we have for any $k$ and $x_{\alpha_k} \in X_m$

\bean &&
(1-\kappa) \|R_{\alpha_{\hat{k}}}T({x_0}-\hat{x}_{\alpha_{\hat{k}}})\|^2\\
&\leq& (1+\kappa)\|R_{\alpha_k}T({x_0}-x_{\alpha_k})\|^2 + 2
{\rm pen}(\alpha_k) + 2
\sup_{k}\{\frac{1}{\kappa}\|R_{\alpha_k}\varepsilon\|^2 - {\rm
pen}(\alpha_k)\}, \eean \vskip .1in

On the other hand, using that is $1\leq R_{\alpha_k}T \leq C$, we
have that for any $x_{\alpha_k} \in X_{m_0}$ and any $k \in \BBn$,

\bean && (1-\kappa) \|{x_0} - \hat{x}_{\alpha_{\hat{k}}}\|^2 \le \,C (1-\kappa) \|{x_0} - x_{\alpha_{k}}\|^2\\
&+& 2\,{\rm pen}(\alpha_k) + 2 \, C_1 \,
\sup_{k}\,\{\|R_{\alpha_k}\varepsilon \|^2 - {\rm
pen}(\alpha_k)\}.\eean As above,  the proof then follows directly
from   lemma \ref{lema1} which characterizes the supremum of the
   empirical process under the linear application as defined by the regularization
family.
\end{proof}

\section{Appendix}\label{slemmas}
 In this section we give some technical lemmas.
The next lemma characterizes the supremum of an empirical process
by the norm of an orthogonal projection.
\begin{lem} \label{lema0}
  \begin{equation}
    \label{proj}
    \sup_{y \in Y_m,\: \| y \|_n=1} |<\varepsilon,y>_n| = \| \Pi_{Y_m}^n \varepsilon \|_n
  \end{equation}

\end{lem}
\begin{proof}
Using the definition of an orthogonal projector, we have
$$ \| \varepsilon - \frac{ 1}{\| \Pi_{Y_m}^n \varepsilon \|_n} \Pi_{Y_m}^n \varepsilon \|^2_n = \min_{ \{ y \in Y_m,\:
  \| y\|_n =1 \}} \| \varepsilon - y \|_n^2.$$
As a consequence we can write:
$$
 \| \varepsilon\|_n^2 - 2 < \varepsilon ,  \frac{ 1}{\| \Pi_{Y_m}^n \varepsilon \|_n} \Pi_{Y_m}^n\varepsilon >_n +  \frac{
   1}{\| \Pi_{Y_m}^n \varepsilon \|_n^2} \| \Pi_{Y_m}^n \varepsilon \|_n^2  =  \min_{ \{y\in Y_m,\:
  \| y\|_n =1 \}} \| \varepsilon\|_n^2 - 2 < \varepsilon, y>_n+ 1 $$
$$2 < \varepsilon -   \Pi_{Y_m}^n \varepsilon ,\frac{ 1}{\| \Pi_{Y_m}^n \varepsilon \|_n} \Pi_{Y_m}^n \varepsilon >_n + 2 < \Pi_{Y_m}^n \varepsilon
,\frac{1}{\| \Pi_{Y_m}^n \varepsilon \|_n} \Pi_{Y_m}^n \varepsilon
>_n  = 2 \sup_{ \{y \in Y_m,\:
  \| y\|_n =1 \}} | < \varepsilon,y>_n| $$
$$ \| \Pi_{Y_m}^n \varepsilon \|_n  =  \sup_{ \{ y \in Y_m,\:
  \| y\|_n =1 \}} | < \varepsilon,y>_n|, $$
which ends the proof.
\end{proof}

  The next result is a  deviation inequality   based on a functional exponential inequality (Theorem 7.4) due to
  \cite{bousquet}
2003. Set $\eta(A)=\sup_{\|u\|=1}\sum_{i=1}^n \varepsilon_i
(A^tu)_i$ for $A:\BBr^n\to \BBr^k$. Let
$$v= \BBe \sum_{i=1}^n\sup_{\|u\|=1} \frac{(A^tu)_i^2}{\rho(A^t
A)} (\frac{\varepsilon_i}{\sigma})^2+2 \BBe \eta(A)/(\sigma
\rho^{1/2}(A^tA)).$$ Then,
\begin{lem}\label{lemabousquet}
$$P(\frac{\eta(A)}{\sigma\rho^{1/2}(A^t A)} > \BBe \frac{\eta(A)}{\sigma \rho^{1/2}(A^t A)}+\sqrt{2 vx}+ x)\le
e^{-x} .
$$
\end{lem}

\begin{proof}   Since the application $u\to A^t u$ is
continuous, we have $\eta(A)= \sup_{u\in S  }\sum_{i=1}^n
\varepsilon_i (A^tu)_i$ for $S  $ some countable subset of the
unit ball. On the other hand,

\bean \sup_{\|u\|=1}[A^tu]_i\le \sup_{\|u\|=1} \|A^t u\| \le
\rho(A). \eean

Thus $\sup_{\|u\|\le 1 }|(A^tu)_i/\rho^{1/2}(A^t A)|\le 1  .$
Also, following the proof of Corollary 5.1 in \cite{Baraud2000}
\bean
&&   \sup_{\|u\|=1}
\frac{(A^tu)_i^2}{\rho(A^t A)}\\
&&\le   \sup_{\|u\|=1}
 \frac{(\sum_{j=1}^m u_j(A^te_j)_i )^2}{\rho(A^t A)}\\
&& \le   \sup_{\|u\|=1}
 \frac{(\sum_{j=1}^m (A^te_j)_i )^2\sum_{j=1}^m u_j^2}{\rho(A^t A)}\\
 && := z_i  .
 \eean

Set $Z=Z(\varepsilon_1,\ldots,\varepsilon_n)=\eta(A)/(\sigma \rho^{1/2}(A^tA))$. Let $\BBe_j$ stand for the conditional expectation
given $\varepsilon_i$ for $i\ne j$.
 Hence, in the proof of Theorem 7.4 in \cite{bousquet} we may bound
 $$|Z-\BBe_jZ|^p\le \frac{|\varepsilon_j|^p}{\sigma^p} \sup_{\|u\|=1}
\frac{(A^tu)_j^2}{\rho(A^t A)}\sup_{\|u\|=1} \max_i
(\frac{(A^tu)_i^2}{\rho(A^t A)})^{p-2}\le (|\varepsilon_j|/\sigma)^p z_j.$$ Thus,
$\BBe |Z-\BBe_j Z|^p\le  z_j p!/2 $. Finally, note that $$\sum_{j=1}^n
z_j=\frac{Tr(A^t A)}{\rho(A^t A)}.$$
   Thus, the proof follows from   Theorem 7.4 in \cite{bousquet}.

\end{proof}

As a corollary, we have the following  lemma
\begin{lem} \label{lema1}
\begin{itemize}
\item There exists a positive constant $d$ that depends on $r/2$ such that the  following inequality holds \bea{desprob} &&P(\eta^2(A)\ge
\sigma^2[Tr(A^tA)+ \rho(A^tA)]r/2(1+L)+\sigma^2 u)\\ \nonumber
&&\le \exp\{-\sqrt{d (1/\rho(A^tA) u+ r/2 L [Tr(A^tA)/\rho(A^t
A)+1] )}\}\eea
\item Set $k_1= d/(\rho(A^tA)\sigma^2)$ and $k_2=d r/2 L
[Tr(A^tA)/\rho(A^tA)+1]$. Then, there exists a constant $C_q$,
which depends only on $q$, such that, \bea{desesp}&& \BBe
[\eta^2(A)-\sigma^2[Tr(A^tA)+ \rho(A^tA)]r/2(1+L) ]_+^q\\
\nonumber &&\le C_q
 k_1^{-q} [ k_2^{q-1/2}+k_2^{q-1}] e^{-\sqrt{  k_2}} \eea  holds.
\end{itemize}
\end{lem}

\begin{proof} As a first step we will bound $v$.    Since $\BBe Z\le \BBe
^{1/2}Z^2$, we have
$$v\le  \BBe  \sum_{i=1}^n z_i (\frac{\varepsilon_i}{\sigma})^2 +2\sqrt{\BBe  \sum_{i=1}^n z_i (\frac{\varepsilon_i}{\sigma})^2 }
\le (1+\nu)\BBe  \sum_{i=1}^n z_i
(\frac{\varepsilon_i}{\sigma})^2 +Tr(A^tA)/\rho(A^tA).$$ Moreover,
following, \cite{Baraud2000} p. 480, for all $p \geq 2$, the following version of Rosenthal's inequality holds:
$$\BBe^{p/2}  \sum_{i=1}^n z_i (\frac{\varepsilon_i}{\sigma})^2 \le  2^{p/2} Tr(A^tA)/\rho(A^tA)\BBe
\frac{|\varepsilon_1|^p}{\sigma^p}.$$ Hence, we have
$$v\le (1+\nu)Tr(A^tA)/\rho(A^tA) +\frac{1}{\nu}$$
and
$$v^2\le 2[ 2^2(1+\nu)^2Tr(A^tA)/\rho(A^tA)\BBe \frac{|\varepsilon^4|}{\sigma^4}
+(\frac{1}{\nu})^2].$$ Set $0<\alpha<1$.
    Choose $\delta$ and $\beta$ such that if
\begin{align*} 2^2 4!\delta^2 (1+1/\alpha)(1-\nu)^2 & <c_1, \\
 2\delta^2(1+1/\alpha) (\frac{1}{\nu})^2 & <c_2
\end{align*}
and $c=\max((1+\beta) \max(c_1,c_2),(1+\beta)(1+\alpha))$,
    then
    $r/2>c$. Let $u>0$ and  without loosing generality,
assume $\sigma=1$.
      Thus, \bean && P(\eta^2(A)\ge
(Tr(A^tA) +\rho(A^tA))r/2(1+L)+u)\\&&\le  P(\eta^2(A)\ge
(Tr(A^tA)(1+\alpha)+(1+1/\alpha)\delta^2
\nu^2\rho(A^t
A))(1+\beta)\\ &&+ [r/2-c] (Tr(A^tA)+\delta^2 v^2\rho(A^t
A))+r/2L (Tr(A^tA) +\rho(A^tA)) +u)\\
&&\le P(\eta^2(A)\ge (Tr(A^tA)(1+\alpha)+(1+1/\alpha)\delta^2
\nu^2
\rho(A^t A)) (1+\beta)\\ &&+
 r/2L (Tr(A^tA) +\rho(A^tA)) +u) \eean
Set $$x'=(1+\frac{1}{\beta})^{-1}\left[\frac{r}{2L} (\frac{Tr(A^tA)}{\rho(A^t A)} +1 )+\frac{u}{\rho(A^tA)} \right].$$ The
last term is equal to
\begin{align*}
&P(\frac{\eta^2(A)}{ \rho(A^t A)}\ge
(\frac{Tr(A^tA)}{\rho(A^t A)}(1+\alpha)\\
& \mspace{20mu}+(1+1/\alpha)v^2\delta^2 )
(1+\beta)  + r/2 L  \frac{Tr(A^tA) }{\rho(A^t A) } +1 )+u)\\
&\mspace{40mu}= P(\frac{\eta^2(A)}{\rho(A^t A)}\ge (\frac{Tr(A^tA)}{\rho(A^t
A)}(1+\alpha)+(1+1/\alpha)v^2\delta^2 )(1+\beta)+
(1+1/\beta)x')
\end{align*}
Finally, we may bound \bean
 &&\le P(\frac{\eta^2(A)}{\rho(A^t A)}\ge (\BBe
\frac{\eta(A)}{\rho^{1/2}(A^tA)}+\delta
v )^2(1+\beta)+(1+1/\beta)x')\\
&&\le P(\frac{\eta^2(A)}{\rho(A^t A)}\ge (\BBe \frac{
\eta(A)}{\rho^{1/2}(A^tA)}+\delta v + \sqrt{x^{\prime }} )^2)
 \\
&&= P(\frac{\eta(A)}{\rho^{1/2}(A^t A)}\ge \BBe
\frac{\eta(A)}{\rho^{1/2}(A^tA)}+ \delta   v  +(1+2/\delta)  x^{\prime\prime})\\
 && \le P(\frac{\eta(A)}{\rho^{1/2}(A^t A)}\ge \BBe
\frac{\eta(A)}{\rho^{1/2}(A^tA)}+ \sqrt{2  v  x^{\prime\prime}} +
x^{\prime\prime})\le
e^{-  x^{''}}\\
&&=e^{-\sqrt{g(A)}}, \eean where we have used
repeatedly that for any constant $c>0, ca^2+1/c b^2\ge 2ab$ and set
$$g(A)=  ((1+1/\beta)^{-1}(1+2/\delta)^2) (r/2 L
[Tr(A^tA)/\rho(A^t A)+1]+u/\rho(A^t A)).$$ Set also
$d= [(1+1/\beta)^{-1}(1+2/\delta)^2  ]^{-1}$ and
$b(A)=Tr(A^tA)/\rho(A^t A)$. Thus we have shown the first part of
the lemma.

Moreover, using the above inequality,

\begin{align*} & {\bf E}[\eta^2(A)-\sigma^2(Tr(A^tA)+ \rho(A^t
A))r/2(1+L)]_+^q \\
& \le \int_0^\infty \sigma^{2q} q u^{q-1}e^{- \sqrt{
  d r/2 L [b(A)+1] + d u/(\rho(A^t A))}} du.
\end{align*}

Consider the change of variable $w=d  u/(\rho(A^t A))+ d r/2 L
[b(A)+1]$, so that
 \begin{align*} & \BBe [\eta^2(A)-\sigma^2(Tr(A^tA)+ \rho(A^t
A))r/2(1+L)]_+^q  \\ \le & \left(\frac{\sigma^2 \rho(A^t
A)}{d}\right)^{q} \int_{ d r/2 L [b(A)+1]}^\infty
 (w- d r/2 L
[b(A)+1])^{q-1} e^{-\sqrt{w }}dw.
\end{align*}

 The last expression is in turn bounded by
 \bean && \left(\frac{\sigma^2 \rho(A^t A)}{d}\right)^{q}
 \int_{d r/2 L [b(A)+1]}^\infty  e^{-\sqrt{w }}[w^{q-1} +
(d r/2 L [b(A)+1])^{q-1}]dw\\
&\le & C_q
 k_1^{-q} [ k_2^{q-1/2}+k_2^{q-1}] e^{-\sqrt{  k_2}}, \eean
ending the proof.

\end{proof}

\thanks
The second author would thank ECOS Nord and Agenda Petr\'oleo de Venezuela for supporting her work.

\bibliographystyle{plain}

\bibliography{essai}

\end{document}